\documentclass[11pt]{article}
\usepackage{amssymb, amsthm, amsmath, amscd}
\newtheorem{thm}{Theorem}[section]
\newtheorem{lem}[thm]{Lemma}
\newtheorem{prop}[thm]{Proposition}
\newtheorem{cor}[thm]{Corollary}
\newtheorem{NN}[thm]{}
\theoremstyle{definition}\newtheorem{df}[thm]{Definition}
\theoremstyle{definition}
\theoremstyle{definition}

\renewcommand{\phi}{\varphi}

\newcommand{\Z}{\mathbb{Z}}
\newcommand{\Q}{\mathbb{Q}}
\newcommand{\R}{\mathbb{R}}

\newcommand{\T}{\mathbb{T}}

\newcommand{\hm}{homomorphism}
\newcommand{\dt}{\delta}
\newcommand{\ep}{\epsilon}
\newcommand{\andeqn}{\,\,\,{\rm and}\,\,\,}
\newcommand{\rforal}{\,\,\,{\rm for\,\,\,all}\,\,\,}
\newcommand{\CA}{$C^*$-algebra}
\newcommand{\SCA}{$C^*$-subalgebra}

\newcommand{\af}{{\alpha}}
\newcommand{\bt}{{\beta}}

\newcommand{\D}{\mathbb D}
\newcommand{\beq}{\begin{eqnarray}}
\newcommand{\eneq}{\end{eqnarray}}
\title{Embedding   Crossed Products into a Unital Simple AF-algebra
}
\author{Huaxin Lin
 }
\date{}

\begin{document}

\maketitle

\begin{abstract}
Let $X$ be a compact metric space and let $\af$ be a homeomorphism on $X.$
Related to a theorem of Pimsner, we show that $C(X)\rtimes_{\af}\Z$ can be embedded into a unital simple AF-algebra if and only if
there is a strictly positive
$\af$-invariant Borel probability measure.

Suppose that $\Lambda$ is a $\Z^d$ action on $X.$ If
$C(X)\rtimes_{\Lambda}\Z$ can be embedded into a unital simple
AF-algebra, then there must exist a strictly positive
$\Lambda$-invariant Borel probability measure. We show that, if in
addition, there is a generator $\af_1$ of $\Lambda$ such that $(X,
\af_1)$ is minimal and unique ergodic, then
$C(X)\rtimes_{\Lambda}\Z^d$ can be embedded into a unital simple
AF-algebra with a unique tracial state.

Let $A$ be a unital separable amenable simple \CA\, with tracial
rank zero and with a unique tracial state which satisfies the
Universal Coefficient Theorem and  let $G$ be a finitely generated discrete abelian group.
Suppose $\Lambda: G\to Aut(A)$ is a
\hm. Then $A\rtimes_{\Lambda} G$ can always be embedded into a
unital simple AF-algebra.

\end{abstract}

\section{Introduction}
Let $X$ be a compact metric space and let $\af$ be a homeomorphism on $X.$
It was proved by Pimsner (\cite{Pi}) that $C(X)\rtimes_{\af}\Z$ can be embedded into an AF-algebra
if and only if $\af$ is pseudo-non-wondering.
Let $A$ be a unital AF-algebra and let $\af\in Aut(A).$ Nate Brown proved the following AF-embedding theorem:
$A\rtimes_\af\Z$ can be embedded into an AF-algebra if and only if $A\rtimes_\af\Z$ is quasidiagonal.
He also gave a $K$-theoretical necessary and sufficient condition for which $A\rtimes_\af\Z$ can be embedded into an AF-algebra. This result has since been improved, at least partially, by Matui (\cite{M1}). For non-$\Z$ actions, Dan Voiculescu
asked when $C(X)\rtimes_{\Lambda}\Z^2$ can be embedded into an AF-algebra.

An AF-algebra may have some infinite feature. Let ${\cal K}$ be the \CA\, of compact operators
on $l^2. $ Then an essential extension $0\to {\cal K}\to E\to M_n\to 0,$ where
$M_n$ is a matrix algebra, gives a unital AF-algebra. This AF-algebra $E$ does not have a faithful tracial state.
In particular, $E$ can not be embedded into a unital simple AF-algebra.
For a unital stably finite 
\CA\, $A$ with a faithful tracial state,  a stronger  and,  perhaps, a more interesting embedding question  is when
$A$ can be embedded into a unital simple AF-algebra. In this note, we will present some progress on
the last question.

We consider  the following question: When can
$C(X)\rtimes_{\Lambda}\Z^d$ be embedded into a unital simple
AF-algebra?  An obvious necessary condition for
$C(X)\rtimes_{\Lambda}\Z^d$ to be embedded into a unital simple
AF-algebra is that there is a strictly positive
$\Lambda$-invariant Borel probability measure on $X$ (see \ref{Dnot} below). We show that
when $d=1,$ if there is such a measure, then, indeed,
$C(X)\rtimes_{\Lambda}\Z^d$ can be embedded into a unital simple
AF-algebra. When $d>1,$  if in addition, there is a generator
$\af_1$ of action $\Lambda$ such that $(X,\af_1)$ is minimal and
unique ergodic, we show that $C(X)\rtimes_{\Lambda}\Z^d$ can be
embedded into a unital simple AF-algebra.

Turn to non-commutative cases, it was shown by N. Brown (\cite{Bn3}) that, among other things,
if $A$ is a UHF-algebra, 
then $A\rtimes_{\Lambda} G$ can be embedded into an AF-algebra, where $G$ is a finitely generated 
discrete abelian group and where 
$\Lambda: \Z^d\to Aut(A)$ is a \hm.  
Let $A$ be a unital  separable
amenable  simple \CA\, with tracial rank zero and with a unique
tracial state which satisfies the UCT and let $\Lambda: G\to
Aut(A).$  We show  $A\rtimes_{\Lambda} G$ can be embedded into a
unital simple AF-algebra. In particular, if $A$ is a unital simple
A$\T$-algebra with a unique tracial state, $A\rtimes_{\Lambda} G$ can
be embedded into a unital simple AF-algebra.

\vspace{0.2in}

{\bf Acknowledgments} This work is partially supported by a grant from National Science Foundation of
U.S.A. The author would like to thank  N. C. Phillips and Hiroki Matui for some helpful conversations.

\section{Preliminaries}

\begin{NN}
{\rm

 Let $A$ be a \CA. Denote by $Aut(A)$ the group of automorphisms on $A.$
 \vspace{0.15in}

 Let $A$ be a stably finite \CA. Denote by $T(A)$ the tracial state space of $A$
 and by $Aff(T(A))$ the normed space of all real affine continuous functions on $T(A).$
 Denote by $\rho_A: K_0(A)\to Aff(T(A))$ the positive \hm\, induced by
 $\rho_A([e])=(\tau\otimes T_k)(e),$ where $e$ is a projection in $A\otimes M_k$ and $T_k$
 is the standard trace on $M_k,$ $k=1,2,....$
 \vspace{0.15in}

 All ideals in this paper are closed two-sided ideals.

}
\end{NN}

\begin{df}\label{Dnot}
{\rm Let $X$ be a compact metric space and let $\af_1,\af_2,...,\af_d$ be homeomorphisms on $X,$
where $d\ge 1$ is an integer.
Denote by $\af_j^*: C(X)\to C(X)$ the automorphism defined by $\af_j^*(f)=f\circ \af_j$ for all $f\in C(X),$
$j=1,2,...,d.$ Suppose that $\af_j\circ \af_i=\af_i\circ \af_j,$ $i,j=1,2,...,d.$ Then
it gives a $\Z^d$ action on $X.$ This gives a \hm, $\Lambda: \Z^d\to Aut(C(X)).$
The crossed product will be denoted by $C(X)\rtimes_\Lambda\Z^d.$ If $d=1,$ we also use
$C(X)\rtimes_{\af_1}\Z$ for the crossed product.
\vspace{0.1in}

 A measure $\mu$ on $X$ is said to be {\it strictly positive $\Lambda$-invariant}  measure,
if $\mu(O)>0$ for any non-empty open subset $O\subset X$ and $\mu(\af_j(E))=\mu(E)$ for
any Borel subset $E\subset X,$
$j=1,2,...,d.$

}
\end{df}

\begin{df}\label{Dz}
{\rm Let $A$ be a unital \CA\, and let $\af\in Aut(A).$ Denote by
$A\rtimes_\af\Z$ the crossed product. In this paper, we will fix a
unitary and denote it by $z_\af$ for which ${\rm ad}\,
z_\af(a)=\af(a)$ for all $a\in A.$ }

\end{df}

\begin{NN}
{\rm Let $X$ be a metric space and let $x\in X.$ For $\dt>0,$ we
will use $O_{\dt}(x)$ for the open ball with center at $x$ and
with radius $\dt.$

For any subset $S\subset X,$ $\overline{S}$ is the closure of $S.$
we denote by $\partial(S)$ the boundary ${\overline{S}}\setminus S$ of $S.$ 
}
\end{NN}

\begin{df}

{\rm Let $X$ be a compact metric space and let $d>0.$
A subset $S\subset X$ is said to be {\it $d$-connected} if for any two pints $x, y\in S,$ there
are $x_0, x_1,...,x_n\in S$ such that
\beq\nonumber
{\rm dist}(x_i, x_{i+1})<d\andeqn x=x_0, y=x_n
\eneq
$i=0,1,...,n-1.$  $S$ is said to be a {\it $d$-connected} component if $S$ is a closed and open subset
of $X$ which is $d$-connected. Clear that every connected component of $X$ is $d$-connected
component of $X.$

}

\end{df}

\begin{df}\label{DU}
{\rm
Denote by ${\cal U}$ throughout this paper the universal UHF-algebra
${\cal U}=\otimes_{n\ge1}M_n.$

Let $\{e_{i,j}^{(n)}\}$ be the canonical matrix units for $M_n.$
Let $u_n\in M_n$ be the unitary matrix such that
${\rm ad}\, u_n(e_{i,i}^{(n)})=e_{i+1,i+1}^{(n)}$ (modulo $n$).
Let $\sigma=\otimes_{n\ge 1} {\rm ad}\, u_n \in Aut({\cal U})$ be the shift
(see for example Example 2.2 of \cite{Bn1}). A fact that we will use in this paper is the following cyclic Rokhlin
property that $\sigma$ has:
For any integer $k>0,$
any $\ep>0$ any finite subset ${\cal F}\subset {\cal U},$  there exist mutually orthogonal projections
$e_1,e_2,...,e_k\in {\cal U}$ such that
\beq
&&\sum_{i=1}^ke_i=1_{\cal U},\\
&& \|xe_i-e_ix\|<\ep\,\rforal x\in {\cal F}\andeqn\\
&&\sigma(e_i)=e_{i+1}, i=1,2,..., k\,\,\, {\rm (}e_{k+1}=e_1{\rm )}.
\eneq
}
\end{df}

We will frequently use the following result.

\begin{thm}{\rm (Theorem 3.4 of \cite{Lngpots})}\label{gpots}

Let $A$ be a unital separable simple \CA\, with tracial rank zero and let $\af\in Aut(A).$
Suppose that $\af$ satisfies the tracial cyclic Rokhlin property.
Suppose also that there is a subgroup $G\subset K_0(A)$ for which $\rho_A(K_0(A))$ is dense
in $Aff(T(A))$ such that $(\af^r)_{*0}|_{G}={\rm id}_G$ for some integer $r\ge 1.$ Then $A\rtimes_\af\Z$ has tracial rank zero.

\end{thm}
If, in addition, $A$ is assumed to be amenable and satisfy the Universal Coefficient Theorem, 
$A\rtimes_\af\Z$ is a unital simple AH-algebra with real rank zero and with no dimension growth,
by the classification theorem of \cite{Lnduke}.

\section{ Embedding into the Cantor systems}

\begin{lem}\label{Lema}
Let $X$ be a compact metric space and $\mu$ be a Borel probability measure on $X.$
Then, for any $x\in X$ and any $\dt>0,$ there exists $0<r<\dt$ such that
$$
\mu(\partial{(O_r(x))})=0.
$$
\end{lem}

This is known. See the proof of 3.2 of \cite{Lndm} for example.

\begin{lem}\label{Lemb}
Let $X$ be an infinite compact metric space and $\Lambda: \Z^d\to Aut(C(X))$ be a  homomorphism.
Suppose that there is a $\Lambda$-invariant strictly positive Borel probability measure $\mu.$ Then, there is a
unital simple AF-algebra $C_0$ with a unique trace $\tau$ and a unital abelian \SCA\, $C_{00}\subset C_0$
 which have the following properties:

{\rm (1)} $\rho_{C_0}(K_0(A))=K_0(A)=\D,$ where $\D$ is a countable divisible dense subgroup of $\R,$

{\rm (2) } there exists a monomorphism $h: C(X)\to C_{00}$ such that
\beq\label{lemb-1}
\tau\circ h(f)=\int_Xf d\mu\,\,\,\rforal f\in C(X),
\eneq

 {\rm (3)} for each $j\in \{1,2,...,d\},$ there exists path of unitaries $\{u_{j,t}: t\in [1,\infty)\}$
in $C_0$  such that
\beq\label{lemb-2}
\lim_{t\to\infty}\|h(f\circ \af_j)-{\rm ad}\, u_{j,t}\circ h(f)\|=0
\eneq
for all $f\in C(X)$ and

{\rm (4)} $C_{00}$ is an AF-algebra with $1_{C_{00}}=1_{C_0}$ and
there exists a sequence of mutually commuting projections $\{e_{n,i}\}$ which generates $C_{00}.$
Moreover
\beq\label{lemb-2+1}
u_{j,n+t}^*e_{n,i}u_{j,n+t}=u_{j,n}^*e_{n,i}u_{j,n} \andeqn
\eneq
\beq\label{lemb-2+2}
u_{j',n+1+t}^*u_{j,n+t}^*e_{n,i}u_{j,n+t}u_{j',t+1}=u_{j,t}^*u_{j',n+1+t}^*e_{n,i}u_{j',n+1+t}u_{j,n+t}
\eneq
for all $n,i,j$ and $t\ge 0.$

\end{lem}

\begin{proof}
Fix a decreasing sequence of positive numbers $\{d_n\}$ for which $\lim_{n\to\infty}d_n=0.$

 Let
$X_{n,1},X_{n,2},...,X_{n,m(n)}$ be disjoint $d_n$-connected components with
$\cup_{i}X_{n,i}=X.$
There are finitely many open subsets $O_1,....,O_{m'(n)}$ such that
$$
\cup_{m=1}^{m'(n)}O_i=X\andeqn \mu(\partial{(O_j)})=0,
j=1,2,...,m'(n).
$$
From this, it is easy to obtain
mutually disjoint Borel subsets $Y'_{n,1},$$ Y'_{n,2},$
$...,Y'_{n,l'(n)}$
 for which each $Y_{n,l'}'$($l'=1,2,...,l'(n)$)  has diameter  less
than $d_n$ and $\cup_{n=1}^{l(n)}Y_{n,l}'=X.$
Moreover, each $Y'_{n,l'}$ has the form $O\cup S,$ where $O$ is an open subset and $S\subset \partial{(O)}$
is a Borel subset and $\mu(\partial(Y'_{n,l'}))=0.$

Furthermore, we assume that $\{Y'_{n,l'}:
l'=1,2,...,l'(n)\}$ is a refinement of $\{ X_{n,1},X_{n,2},...,X_{n,m(n)}\}.$ Let
$$
Y_{n,1}, Y_{n,2},..., Y_{n,l(n)}
$$
be mutually disjoint Borel subsets for which each $Y_{n,l}$ ($l=1,2,...,l(n)$) has diameter  less than $d_n$ and the collection
of finite union of these $Y_{n,l}$'s  contains
$$
\{X_{n,i},\,\af_j^{-1}(X_{n,i}): i=1,2,...,m(n)\}\cup \{ Y'_{n,l}, \af_j^{-1}(Y'_{n,l}): l=1,2,...,l(n)\}.
$$
$j=1,2,...,d,$ $n=1,2,....$
Moreover, 
We may assume that we also assume that each $Y_{n,l}$ has the form $O\cup S,$ where $O$ is an open and $S\subset \partial{(O)}$
is a Borel subset as well as $\mu(\partial(Y_{n,l}))=0.$

We may assume that 
$$
Y'_{n+1,1},Y'_{n+1, 2},...,Y'_{n+1,l'(n+1)}
$$
is a refinement of $\{Y_{n,1}, Y_{n,2},..., Y_{n,l(n)}\}.$
By induction, 
it follows that
 the partition
$\{Y_{n+1,1}, Y_{n+1,2},..., Y_{n+1,l(n+1)}\}$ is a refinement of $\{Y_{n,1}, Y_{n,2},..., Y_{n,l(n)}\},$
$n=1,2,....$
Define $\D$ to be a countable divisible group of $\R$ which contains
$\Q$ and
$\{\mu(Y_{n,i}): i=1,2,...,l(n), n=1,2....\}.$ Let $C_0$ be a unital separable simple AF-algebra
with unique trace $\tau$ such that $(K_0(C_0), K_0(C_0)_+, [1_{C_0}])=(\D,\D_+.1).$

There are mutually orthogonal projections $\{e_{n,i}: i=1,2,...,l(n)\}$ such that
$\tau(e_{n,i})=\mu(Y_{n,i}),$ $i=1,2,...,l(n).$
It follows that $\sum_{i}e_{n,i}=1_{C_0}.$ 
 For a fixed $n,$ if $Y_{n+1,j}\subset Y_{n,i},$
one can obtain projections
$\{e_{n+1,j}\}$ in  the \SCA\, $e_{n,i}C_0e_{n,i}.$ Moreover, if for some finite subset $J,$
$\cup_{j\in J} Y_{n+1,j}=Y_{n,i},$
then $\sum_{in J}e_{n+1,j}=e_{n,i}.$ Therefore, we may assume that $\{e_{n,i}: i=1,2,...,l(n),
n=1,2,...\}$ is also a set of pairwise commuting projections.  Let $C_{00}$ be the commutative \SCA\,
of
$C_0$ generated by
$\{e_{n,i}: i=1,2,...,l(n), n=1,2,...\}.$
Define $h_n: C(X)\to C_{00}$ by
\beq\label{Fh}
h_n(f)=\sum_{i=1}^{l(n)} f(x_{n,i}) e_{n,i}\,\rforal\, f\in C(X),
\eneq
where $x_{n,i}\in Y_{n,i},$ $i=1,2,...,l(n)$ and $n=1,2,....$
For any $m>n$ and $i\in \{1,2,...,l(n)\},$ there exists $J(i)\subset \{1,2,...,l(m)\}$
such that $\sum_{j\in J(i)}e_{m,j}=e_{n,i}.$ We have
\beq\label{Fh2}
h_n(f)-h_m(f)=\sum_{i=1}^{l(n)} \sum_{j\in J(i)}(f(x_{n,i})-f(x_{m,j}))e_{m,j}\,\,\,\rforal\, f\in C(X).
\eneq
Note that $x_{n,i}, x_{m,j}\in Y_{n,i}.$ Since $\lim_{n\to\infty}d_n=0,$ we conclude that,
for each $f\in C(X),$ $\{h_n(f)\}$ is Cauchy in $C_{00}.$
Define $h(f)=\lim_{n\to\infty} h_n(f)$ for all $f\in C(X).$

Write $C_{00}=C(Y),$ where $Y$ is a compact totally disconnected metric space.
Let $\Omega_{n,i}$ be the clopen subset corresponding to $e_{n,i}.$ 
Let $f_1\in C(X)$ such that $f_1(x)=0$ for all $x\not\in Y_{n,i}^o,$ the interior of $Y_{n,i}$ and let $f_2\in C(X)$ such that 
$f_2(x)=1$ for all $x\in \overline(Y_{n,i}).$  It follows from (\ref{Fh}) and (\ref{Fh2}) that
\beq\label{Fh3}
h(f_1)e_{n,i}=h(f_1)\andeqn h(f_2)e_{n,i}=e_{n,i}
\eneq
which will be used later.

Fix $n\ge 2.$ For each $j$ and $i,$ there is a finite subset $J(i,n,j)\subset \{1,2,...,l(n+1)\}$ and a
projection with the form $\sum_{i'\in J(n,i,j)}e_{n+1,i'}$ which corresponds to $\af_j^{-1}(Y_{n,i}).$
Denote this projection by $p_{n,i,j}.$
There is $\xi_{n,i}\in \af^{-1}(Y_{n,i})$ such that $\af_j(\xi_{n,i})=x_{n,i}.$
Since $\lim_{n\to\infty}d_n=0,$ we have
\beq\label{NL-1}
\lim_{n\to\infty} \|\sum_{k=1}^{l(n+1)}f(\af_j(x_{n+1,i}))e_{n+1,k}-\sum_{i=1}^{l(n)}f(x_{n,i})p_{n,i,j}\|=0
\eneq
for all $f\in C(X).$ Fix $j$ and fix $n\ge 2,$ there is a unitary $u_{j,n}\in C_0$ such that
\beq
u_{j,n}^*e_{n,i}u_{j,n}=p_{n,i,j}, j=1,2,...,d
\eneq
$i=1,2,...,l(n).$

Suppose that $u_{j,n}$ has been defined.
Suppose that
$S(i,n)\subset  \{1,2,...,l(n+1)\}$ such that
$\sum_{k\in S(i,n)} \af_j^{-1}(Y_{n+1,k})=\af_j^{-1}(Y_{n,i}).$
Then, for $k\in S(i,n),$
\beq
u_{j,n}^*e_{n+1,k}u_{j,n}\le p_{n,i,j} \andeqn \sum_{k\in S(i,n)}u_{j,n}^*e_{n+1,k}u_{j,n}=p_{n,i,j}.
\eneq

Let $J(k,n+1)\subset \{1,2,...,l(n+2)\}$ so that
$\cup_{i'\in J(k,n+1)}Y_{n+2,i'}=\af_j^{-1}(Y_{n+1,k}).$ Let
\beq
p_{n+1, k,j}=\sum_{i'\in j(k,n+1)}e_{n+2,i'}.
\eneq
We have
\beq
\tau(u_{j,n}^*e_{n+1,k}u_{j,n})=\tau(p_{n+1,k,j}).
\eneq
Therefore there  is a partial isometry $w(j,n+1,i,k)\in p_{n,i,j}C_0p_{n,i,j}$ such that
($k\in S(i,n)$)
\beq\nonumber
w(j,n+1,i,k)^*w(j,n+1,i,k)&=&p_{n+1,k,j}\,\,\,\,\,\,\,\andeqn\\
\ w(j,n+1,i,k)w(j,n+1,i,k)^*&=&u_{j,n}^*e_{n+1,k}u_{j,n}.
\eneq
Put $w(j,n+1,i)=\sum_{k\in S(i,n)}w(j,n+1,i,k).$
Since
\beq
\sum_{k\in S(i,n)} p_{n+1,k,j}=p_{n,i,j},
\eneq
$w(j,n+1, i)$ is  a unitary in $p_{n,i,j}C_0p_{n,i,j}.$
There is a path of unitaries in  $\{W(j,n+1,t,i):t\in [0,1]\}$ in $p_{n,i,j}C_0p_{n,i,j}$ with length no more than $\pi$ for which
\beq
W(j,n+1,i,0)=p_{n,i,j}\andeqn W(j, n+1,i,1)=w(j,n+1,i).
\eneq
Define
\beq
U_{j,n}(t)=\sum_{i=1}^{l(n)}W(j,n+1,i,t)u_{j,n}\andeqn u_{j,n+1}=\sum_{i=1}^{l(n)}W(j,n+1,1)u_{j,n}.
\eneq
 Then
 \beq
&&U_{j,n}(0)=u_{j,n}, \,\,\, U_{j,n}(1)=u_{j,n+1},\\\label{03+3}
&& u_{j,n+1}^*e_{n+1,k}u_{j,n+1}=p_{n+1,k,j}\andeqn U_{j,n}(t)^*e_{n,i}U_{j,n}(t)=p_{n,i,j}
\eneq
Define $u_{j,s}=U_{j,n}(s-n)$ for $s\in [n,n+1).$
Thus, for any $n$ and $i,$ 
\beq\label{Fh4}
u_{j,s}^*e_{n,i}u_{j,s}=u_{j,n}^*e_{n,i}u_{j,n}
\eneq
for any $s\ge n.$

Suppose that $J'(i, n+1, j,j_1)\subset  \{1,2,...,l(n+2)\}$ such that
$\sum_{i_1\in J'(i,n+1, j,j_1)} e_{n+2,i_1}$ corresponds to $\af_{j_1}^{-1}\circ \af^{-1}_j(Y_{n,i}).$
Then
since $\af^{-1}_{j'}\circ \af^{-1}_j(Y_{n,i})=\af^{-1}_j\circ \af^{-1}_{j'}(Y_{n,i}),$
\beq
\sum_{i_1\in J'(i,n+1, j,j_1)} e_{n+2,i_1}=\sum_{i_2\in J'(i,n+1, j_1,j)} e_{n+2,i_2}
\eneq
This implies that
\beq\label{03+1}
u_{j_1,n+1}^*u_{j,n}^*e_{n,i}u_{j,n}u_{j_1,n+1}=u_{j,n}^*u_{j_1, n+1}^*e_{n,i}u_{j_1,n+1}u_{j,n}.
\eneq
By applying (\ref{NL-1}), we obtain that
\beq
\lim_{t\to\infty}\|h(f\circ \af_j)-{\rm ad}u_{j,t}\circ h(f)\|=0
\eneq
for all $f\in C(X)$ and $j=1,2,...,d.$ We also have
\beq
\tau\circ h(f)=\int_Xfd\mu
\eneq
for all $f\in C(X).$
From (\ref{03+1}) and (\ref{03+3}),  one also has  both (\ref{lemb-2+1}) and (\ref{lemb-2+2}).

\end{proof}

\begin{cor}\label{LembC}
Let $X$ be a compact metric space and let $\Lambda: \Z^d\to Aut(C(X))$ be a \hm. Suppose that there is
a strictly positive $\Lambda$-invariant Borel probability measure $\mu.$ Then there is an embedding
$h_1: C(X)\to C(Y),$ where $Y$ is a compact totally disconnected metric space and there is an embedding $h_2: C(Y)\to C_0,$ where
$C_0$ is a unital simple AF-algebra with a unique tracial state $\tau$ satisfying the following:

{\rm (1)} For each $j=1,2,...,d,$ there is a path of unitaries $\{u_{j,t}: t\in [1, \infty)\}$ in
$C_0$
 such that
 \beq
 \lim_{t\to \infty}\|{\rm ad}\, u_{j,t}\circ h_2\circ h_1(f)-\af_j^*\circ h_1(f)\|=0
 \eneq
 for all $f\in C(X).$

 {\rm (2)}
 \beq
 \lim_{t\to\infty} {\rm ad}\, u_{j,t}(h_2(g))=\bt_j^*(g)
 \eneq defines an automorphism on $C(Y).$

 Moreover $\bt_1^*, \bt_2^*,...,\bt_d^*$ defines a $\Z^d$ action ${\overline{\Lambda}}$ on $C(Y)$ and
$$
\bt_j^*\circ h_1=h_1\circ \af_j^*,\,\,\,
j=1,2,...,d.
$$

\end{cor}

\begin{proof}
We will use the notation in the proof of \ref{Lemb}.  Since $C_0$ is a  unital
separable commutative AF-algebra, there is a totally disconnected compact metric space $Y$ such that
$C_{00}\cong C(Y).$ Note that (1) directly follows from \ref{Lemb}.
To see {\rm (2)}, we note that, for each projection $e\in C_{00},$ as in the proof of \ref{Lemb},
by (\ref{Fh4}),
there is $t_0\ge 1$ such that
\beq
u_{j,t}^*eu_{j,t}=u_{j,s}^*eu_{j,s}
\eneq
for all $t,s\ge t_0$ and $j=1,2,...,d.$ In particular,  $u_{j,t}^*eu_{j,t}$ are in $C_{00}$ for $t\ge
t_0.$  It follows that
$\{u_{j,t}^*gu_{j,t}\}$ converges to an element in $C_{00}.$
Thus $\lim_{t\to\infty} {\rm ad}\, u_{j,t}(h_2(g))=\bt_j^*(g)$
 defines an automorphism on $C(Y).$
 It follows from (4) of \ref{Lemb} that $\bt_1^*,\bt_2^*,...,\bt_d^*$ generate a $\Z^d$ action
 $\overline{\Lambda}$ on $C(Y).$
  Since
 \beq
 \bt_j^*(h_1(f))=\lim_{t\to\infty} {\rm ad}\, u_{j,t}(h_2\circ h_1(f))= h_1\circ \af_j^*(f)
 \eneq
 for all $f\in C(X),$
  $\bt_j\circ h_1=h_1\circ \af_j^*,$ $j=1,2,...,d.$ So
 the last part of the corollary follows.

\end{proof}

\begin{cor}\label{C03}
Let $X$ be a compact metric space and let $\Lambda: \Z^d\to
Aut(C(X))$ be a \hm. Suppose that there is a strictly positive
$\Lambda$-invariant Borel probability measure $\mu.$ Let $C_0,$
$C(Y)\cong C_{00},$ $h_1: C(X)\to C(Y)$ and $h_2: C(Y)\to C_0$ be
as in \ref{LembC} and as constructed in \ref{Lemb}. Denote by
$\bt_j: Y\to Y$ the homeomorphism induced by the automorphism
$\bt_j^*.$ Suppose that $\af_1$ is minimal. Then $\bt_1: Y\to Y$
is also minimal and $Y$ is homeomorphic to the Cantor set, if $X$ is infinite.  If
$(X, \af_1)$ has a unique $\af_1$-invariant Borel probability
measure, then $(Y, \bt_1)$ has a unique $\bt_1$-invariant Borel
probability measure.
\end{cor}

\begin{proof}
Fix a continuous surjective map $s: Y\to X$ for which
$h_1(f)(y)=f\circ s (y)$ for all $y\in Y.$ Fix $\xi\in Y$ and let
$\zeta=s(\xi).$ We will keep notations used in the proof of
\ref{Lemb}.

First we claim the following: If $\dt>0$ and if $
\overline{Y_{n,i}}\cap O_{\dt}(\zeta)=\emptyset,$ then there is no
$y\in \Omega_{n,i}$ such that $s(y)=\zeta,$ where
$\Omega_{n,i}\subset Y$ is the clopen set associated with
$e_{n,i}.$

Let $O_1, O_2\subset X$ be open subsets containing
${\overline{Y_{n,i}}}$ such that ${\overline{O_1}}\subset O_2$ and
$O_2\cap O_{\dt/2}(\zeta)=\emptyset.$ Let $f_1\in C(X)$ such that
$0\le f_1(x)\le 1$ for all $x\in X,$ $f_1(x)=1$ if $x\in O_1$ and
$f_1(x)=0$ if $x\not\in O_2.$ Let $f_2\in C(X)$ such that $0\le
f_2(x)\le 1$ for all $x\in X,$ $f_2(x)=1$ if $x\in
O_{\dt/4}(\zeta)$ and $f_2(x)=0$ if $x\not\in O_{\dt/2}(\zeta).$
Then $f_1f_2=0.$ So $h_1(f_1)h_1(f_2)=0.$ Note that, by (\ref{Fh3}),
$h_1(f_1)e_{n,i}=e_{n,i}.$ So $h_1(f_2) e_{n,i}=0.$ If $y\in
\Omega_{n,i}$ such that $s(y)=\zeta.$ Then
$h_1(f_2)(y)e_{n,i}(y)=f_2(s(y))e_{n,i}(y)\not=0.$ This is a
contradiction. The claim is proved.

Now suppose that $\af_1$ is minimal.  To show that $\bt_1$ is
minimal, let $\Omega_{n,i}$ be a clopen subset associated with the
projection $e_{n,i}.$ We will show that there exists an integer
$N$ such that $\bt_1^N(\xi)\in \Omega_{n,i}.$ Let $Y_{n,i}^o$ be
the interior of $Y_{n,i}.$ Choose $x\in Y_{n,i}^o$ and $\ep>0$
such that $\overline{O_{\ep}(x)}\subset Y_{n,i}^o.$ Since $\af_1$
is minimal, there is an integer $N$ such that $\af_1^N(\zeta)\in
O_{\ep/4}(x).$ There is $\dt>0$ such that
 \beq\label{C03+0}
 \af_1^N(O_{\dt}(\zeta))\subset O_{\ep/4}(x).
 \eneq

Let $d_n$ be as in the proof of \ref{Lemb}. Choose $m$ so that
$d_m<\dt/32.$ Note that, for some $i',$ $\xi\in \Omega_{m,i'}.$ By
the claim,
$$
{\overline{Y_{m,i'}}}\cap O_{\dt/32}(\zeta)\not=\emptyset.
$$
Therefore, since the diameter of $Y_{m,i'}$ is no more than $d_m,$
\beq\label{C03+1}
 {\overline{Y_{m,i'}}}\subset
O_{\dt/8}(\zeta). \eneq
 Let $g\in C(X)$ such that $0\le g(x)\le 1$
for all $t\in X,$ $g(t)=1$ if $t\in {\overline{Y_{m,i'}}}$ and
$g(t)=0$ if $t\not\in O_{\dt/8}(\zeta).$ Let $g_1\in C(X)$ such
that $0\le g_1(t)\le 1,$ $g_1(t)=1$ if $t\in O_{\ep/4}(x)$ and
$g_1(t)=0$ if $t\not\in O_{\ep/2}(x).$

 Then (since $Y_{n,i}^o\supset \overline{O_{\ep}(x)}$), by (\ref{Fh3}),
$$
h_1(g)e_{m,i'}=e_{m,i'} \andeqn h_1(g_1)e_{n,i}=h_1(g_1).
$$
Therefore (since $\bt_1^*\circ h_1=h_1\circ \af_1^*$)
 \beq\label{C03+3}
h_1\circ(\af_1^*)^N(g)(\bt_1^*)^N(e_{m,i'})=(\bt_1^*)^N(e_{m,i'}).
\eneq
 By (\ref{C03+1}) and (\ref{C03+0}),
\beq\label{C03+4}
 h_1\circ(\af_1^*)^N(g)
h_1(g_1)=h_1\circ(\af_1^*)^N(g).
\eneq
We also have, by applying
(\ref{C03+3}) and (\ref{C03+4}),

\beq\nonumber
e_{n,i}(\bt_1^*)^N(e_{m,i'})&=&e_{n,i}h_1\circ(\af_1^*)^N(g)(\bt_1^*)^N(e_{m,i'})\\\nonumber
&=&e_{n,i}h_1(g_1)h_1\circ(\af_1^*)^N(g)(\bt_1^*)^N(e_{m,i'})\\\nonumber
&=&h_1(g_1)h_1\circ(\af_1^*)^N(g)(\bt_1^*)^N(e_{m,i'})\\\nonumber
&=&h_1\circ(\af_1^*)^N(g)(\bt_1^*)^N(e_{m,i'})\\\nonumber
&=&(\bt_1^*)^N(e_{m,i'})
 \eneq
This implies that
$$
\bt_1^N(\Omega_{m,i'})\subset \Omega_{n,i}.
$$
Hence
$$
\bt_1^N(\xi)\in \Omega_{n,i}.
$$
Therefore $\bt_1$ is minimal.

Since $(Y, \bt_1)$ is minimal, $Y$ can not have isolated points.
It follows that $Y$ is an infinite  compact totally disconnected
perfect metric space. Therefore $Y$ is homeomorphic to a Cantor
set.

Now suppose that $(X, \af_1)$  has a uniquely $\af_1$-invariant Borel probability measure.
Note that $\mu$ induced by $\tau\circ h_2$ is a $\bt_1$-invariant Borel probability measure on $Y.$
Let $\mu_1$ be another $\bt_1$-invariant Borel probability measure on $Y.$ Let $\tau_1$ be the tracial state
of $C(Y)$ defined by
$$
\tau_1(g)=\int_X gd \mu_1\,\,\,{\rm for}\,\,\, g\in C(Y).
$$
Then $\tau_1\circ h_1$ gives an $\af_1^*$-invariant tracial state on $C(X).$
It follows that
$$
\tau_1\circ h_1(f)=\int_X f d\mu \,\,\,{\rm for}\,\,\,f\in C(X).
$$
For any $Y_{n,i},$ by (\ref{Fh3}),
\beq\label{C03C+100}
\tau_1(e_{n,i})\ge \sup\{\int_X f d\mu : f\in C(X), 0\le f(x)\le 1, f(x)=0\,\,\,{\rm if}\,\,\, x\not\in Y_{n,i}^o\}=\mu(Y_{n,i}^o).
\eneq
Similarly,
\beq\label{C03C+101}\nonumber
\tau_1(1-e_{n,i}) &\ge &\sup\{\int_X f d\mu : f\in C(X), 0\le f(x)\le 1, f(x)=0\,\,\,{\rm if}\,\,\, x\in \overline{Y_{n,i}}\}\\
&=&\mu(X\setminus \overline{Y_{n,i}}).
\eneq
Note, since $\mu(\partial(Y_{n,i}))=0,$
$$
\mu(X\setminus \overline{Y_{n,i}})=\mu(X\setminus Y_{n,i}) \andeqn \mu(Y_{n,i})=\mu(Y_{n,i}^o).
$$
We also have $\tau_1(e_{n,i})+\tau_1(1-e_{n,i})=1$ and
$\mu(Y_{n,i})+\mu(X\setminus Y_{n,i})=1.$ Thus, by (\ref{C03C+100}) and (\ref{C03C+101}),
$$
\tau_1(e_{n,i})=\mu(Y_{n,i}).
$$
We also have $\tau(e_{n,i})=\mu(Y_{n,i}).$
Since $\{e_{n,i}\}$ generates $C(Y),$ we conclude that
$$
\tau_1=\tau\circ h_2.
$$
Thus $Y$ has only one $\bt_1$-invariant Borel probability measure.

\end{proof}

\section{Unital simple AF-embedding}


First we would like to point out that, in general, a \CA\, which
can be embedded into a unital  AF-algebra may not be embedded into
a unital simple AF-algebra. A simplest example is to consider a
unital AF-algebra which is defined by  a unital essential
extension:
\beq\nonumber
0\to {\cal K}\to E\to M_n\to 0
\eneq
$E$ can not be embedded into any unital simple AF-algebra since it does not admit a faithful
tracial state.

\begin{lem}\label{PP}{\rm (Pasnicu-Phillips)}

Let $A$ be a unital \CA\, and let $\af\in Aut(A)$ be an
automorphism which has the following version of cyclic Rokhlin property:
for any $\ep>0,$ any finite subset ${\cal F}\subset A$ and any
integer $n>1,$ there exist mutually orthogonal projections
$p_1,p_2,,....,p_n$ in $A$ such that

{\rm (1)} $\|ap_i-p_ia\|<\ep$ for all $a\in {\cal F}$ and
$i=1,2,...,n,$

{\rm (2)} $\|\af(p_i)-p_{i+1}\|<\ep$ for $i=1,2,...,n-1$ and

{\rm (3)} $\sum_{i=1}^np_i=1.$

Then, for any $a\in A\rtimes_\af\Z \setminus \{0\}$ and any
$\eta>0,$ there are mutually orthogonal projections
$e_1,...,e_m\in A$ such that
\beq\label{PP-2}
\sum_{i=1}^me_m=1\,\,\, and\,\,\,
\|E(a)-\sum_{i=1}^me_iae_i\|<\eta
\eneq
where $E: A\rtimes_\af\Z\to A$ is the canonical conditional
expectation. Consequently, if $I\subset A\rtimes_\af\Z$ is a
proper ideal, then $I\cap A$ is a proper closed $\af$-invariant
ideal of $A.$
\end{lem}

\begin{cor}\label{PP1}
In the situation in \ref{PP}, if $\tau$ is a tracial state of
$A\rtimes_\af\Z,$ then there exists a $\af$-invariant tracial
state $\tau_1$ of $A$ such that
\beq
\tau(b)=\tau_1\circ E(b)\,\rforal b\in A\rtimes_\af\Z,
\eneq
where $E: A\rtimes_\af\Z\to A$ is the standard conditional
expectation.
\end{cor}

\begin{proof}
Let $\ep>0.$ By \ref{PP}, there are mutually orthogonal
projections $e_1,...,e_m\in A$ such that
\beq\label{PP-1}
\sum_{i=1}^me_m=1\,\,\, and\,\,\,
\|E(a)-\sum_{i=1}^me_iae_i\|<\ep
\eneq
Let $\tau$ be a tracial state on $A\rtimes_\af\Z.$ Then
\beq
\tau(e_iae_i)=\tau(ae_i^2)=\tau(ae_i)\,\,\,i=1,2,...,m.
\eneq
Therefore
\beq
\tau(\sum_{i=1}^me_iae_i)&=&\sum_{i=1}^m \tau(e_iae_i)\\
&=&\sum_{i=1}^m\tau(ae_i)=\tau(\sum_{i=1}^mae_i)=\tau(a)
\eneq
Let $\tau_1=\tau|_A,$ where we identify $A$ with its image of
natural embedding in $A\rtimes_\af\Z.$ Then
\beq
|\tau(a)-\tau(E(a))|=|\tau(\sum_{i=1}^me_iae_i)-\tau_1\circ
E(a)|<\ep
\eneq
for all $a\in A\rtimes_\af\Z.$ It follows that
\beq
\tau=\tau_1\circ E.
\eneq

\end{proof}

\begin{thm}\label{EmbM}
Let $X$ be a compact metric space and let $\af$ be a homeomorphism. Then the following are
equivalent.

{\rm (1)} $C(X)\rtimes_\af
\Z$ can be embedded into a unital simple AF-algebra.

{\rm (2)} There is a strictly positive $\af$-invariant Borel probability measure on
$X.$
\end{thm}

\begin{proof}
Suppose that there is a monomorphism $\phi: C(X)\rtimes_\af\Z\to C,$ where $C$ is a unital simple
AF-algebra. Let $\tau$ be a tracial state of $C.$ Then $\tau\circ \phi$ gives a strictly positive
$\af$-invariant Borel probability measure. Thus (1) implies (2).

Now we apply \ref{LembC} and a result of N. Brown to prove (2) $\Rightarrow$ (1). It follows from the of
\ref{LembC} that there is a monomorphism
$j_1: C(X)\to C(Y),$ where is
$Y$ is a compact totally disconnected space
 and a homeomorphism $\bt$ on $Y$ such that
$$
\bt\circ j_1=j_1\circ \af.
$$
Note
 that $C(Y)=C_{00}$ is an AF-algebra.
Let $C_0$ be as in \ref{Lemb}. Let ${\cal U}$ and $\sigma$ be as
in \ref{DU}. Denote by $j_2: C(Y)\to C_0\otimes {\cal U}$  the
composition of the embeddings. Then
\beq
(({\rm
id}_{C_0}\otimes \sigma)\circ j_2)_*=(j_2\circ \bt)_*.
\eneq
Note that ${\rm id}_{C_0}\otimes \sigma$ has the cyclic Rokhlin property of \ref{DU}.
 It follows from 2.8 of \cite{Bn1} that there
is a unitary
$v\in C(Y),$ a unitary
$u\in  C_0\otimes{\cal U}$ and a monomorphism
$\phi: C(Y)\to C_0\otimes {\cal U}$ such that
\beq
({\rm ad}\, u\otimes \sigma)\circ \phi=\phi\circ {\rm ad}\, v\circ
\bt={\rm ad}\, \phi(v)\circ \phi\circ \bt.
\eneq
Let $\gamma: C_0\otimes{\cal U} \to C_0\otimes {\cal U}$ be
defined by ${\rm ad}\, (\phi(v^*))\circ( u\otimes \sigma).$ Then
$\gamma\circ \phi=\phi\circ \bt.$ On $C_0\otimes {\cal U}\otimes
{\cal U}$ define  ${\widetilde \gamma}=\gamma \otimes \sigma.$ We
obtain an injective \hm\,
\beq
{\tilde \phi}: (C(Y)\otimes {\cal U})_{\bt\otimes \sigma}\Z\to (C_0\otimes {\cal
U}\otimes {\cal U})\rtimes_{\widetilde
\gamma}\Z.
\eneq
Note that $({\tilde \gamma})_{*0}={\rm id}_{K_0(C_0\otimes {\cal U}\otimes {\cal U})}.$  Since $\widetilde \gamma$ has
the cyclic Rokhlin property in
\ref{DU}, by
\ref{gpots}, $C_0\otimes{\cal U}\otimes {\cal U}\rtimes_{\widetilde \gamma}\Z$ has tracial rank zero.
Since it also satisfy the UCT, it is a unital simple AH-algebra with real rank zero and with no
dimension growth. Since $C_0$ has a unique tracial state, so does $C_0\otimes {\cal U}\otimes {\cal U}.$ It follows
from
\ref{PP1} that
$C_0\otimes {\cal U}\otimes {\cal U}\rtimes_{\widetilde \gamma}\Z$ has a unique tracial state.
Let $\D$ be the tracial range of
$K_0(C_0\otimes {\cal U}\otimes {\cal U}\rtimes_{\widetilde \gamma}\Z).$ Let $D$ be the
unital simple AH-algebra with $(K_0(D), K_0(D)_+, [1_D])=(\D, \D_+, 1).$ It follows that there is an
injective
\hm\, $\phi_1: C_0\otimes {\cal U}\otimes {\cal U}\rtimes_{\widetilde \gamma}\Z\to D.$ We then first  embed
$C(X)\rtimes_\af\Z$ into $C(Y)\otimes{\cal U}\rtimes_{\bt\otimes \sigma}\Z$ and embed the latter into
$C_0\otimes {\cal U}\otimes {\cal U}\rtimes_{\widetilde \gamma}\Z.$ By composing this embedding with $\phi_1,$ we obtain the
desired embedding.
\end{proof}

\begin{prop}\label{2T2}
Let $A$ be a unital separable \CA. Then the following are equivalent

{\rm (1)} $A$ can be embedded into a unital simple AF-algebra;

{\rm (2)} $A$ can be embedded into a unital simple AF-algebra with a unique tracial state;

{\rm (3)} $A$ can be embedded into a unital separable  amenable simple \CA\, with tracial rank zero  which satisfies
the UCT.

\end{prop}

\begin{proof}
Suppose (1) holds. Let $B$ be a unital simple AF-algebra and let $h: A\to B$ be an embedding. By
replacing
$B$ by $h(1_A)Bh(1_A),$ we may assume that $h$ is unital.
Let $\tau\in T(B)$ be a tracial state. Define $r: K_0(B)\to \R$ by $r([e])=\tau([e]).$ Let
\beq
\D=\{r(x): x\in K_0(B)\}.
\eneq
Then it is known that $\D$ must be a countable dense subgroup of $\R.$ Let $D$ be a unital simple
AF-algebra such that $(K_0(D), K_0(D)_+, [1_D])=(\D, \D_+, 1).$ There is a unital monomorphism
$h_1: B\to D.$ Thus $h_1\circ h: A\to D$ gives an embedding of $A$ into a unital simple AF-algebra
with a unique tracial state.

That (2) $\Rightarrow$ (3) is obvious.

Suppose that (3) holds. Let $B_1$ be a a unital separable amenable simple \CA\, with tracial rank zero
which satisfies the UCT and let $\phi: A\to B_1$ is an embedding.
Since $B_1\otimes {\cal U}$ also satisfies the UCT, by the classification theorem (\cite{Lnduke}),
there is a monomorphism $\phi_2: B_1\to C$ for some unital simple AF-algebra. This
implies that (1) holds.

\end{proof}

It should be noted that, when $d=1,$  Matui (\cite{M1}) proved
that $A$ can be embedded into an AF-algebra.

\begin{thm}\label{TL}
Let $A$ be a unital separable simple amenable \CA\, with tracial
rank zero and with a unique tracial state which satisfies the UCT.
Suppose that $\Lambda: \Z^d\to Aut(A)$ is a \hm. Then
$A\rtimes_{\Lambda}\Z^d$ can be embedded into a unital simple
AF-algebra
\end{thm}

\begin{proof}
Suppose that $\Lambda$ is determined by $d$ mutually commuting
automorphisms $\af_1, \af_2,...,\af_d.$ Let ${\cal U}$ and
$\sigma$ be as in \ref{DU}. Put $A_1=A\otimes{\cal U}$ and
$\gamma_1=\af_1\otimes \sigma.$ Define
$B_1=A_1\rtimes_{\gamma_1}\Z.$ Then $\gamma_1$ has the cyclic
Rokhlin property (\ref{DU}). For any projection with the form
$1_A\otimes e,$ where $e\in {\cal U}$ is a projection, one has
$$
[\gamma_1(1_A\otimes e)]=[1_A\otimes e]\,\,\, {\rm in}\,\,\,
K_0(A_1).
$$
Let $G_1\subset K_0(A_1)$ be the subgroup generated by projections
of the form $1_A\otimes e.$ Then
$$
(\gamma_1)_{*0}|_{G_1}={\rm id}_{G_1}.
$$
Let $\tau$ be the unique tracial state on $A_1.$  Then
$$
\{\tau(1_A\otimes e): e \,\,\,{\rm projections\,\,\,in}\,\,\,
M_k({\cal U}), k=1,2,...\}
$$
is dense in $\R.$  Since $A$ has a unique tracial state, so does
$A_1.$ It follows that $\rho_{A_1}(G_1)$ is dense in
$Aff(T(A_1)).$ Then, by \ref{gpots}, $A_1\rtimes_{\gamma_1}\Z$ has
tracial rank zero and satisfies the UCT. We also note that there
is an embedding:
\beq\label{F1}
A\rtimes_{\af_1}\Z\to (A\otimes {\cal
U})\rtimes_{\gamma_1}\Z=A_1\rtimes \Z.
\eneq
In particular, if
$d=1,$ by \ref{2T2}, the theorem follows.

If $d>1,$ put $B_1=A_1\otimes_{\gamma_1}\Z$ and $A_2=B_1\otimes
{\cal U}.$ It follows from \ref{PP1}, since $\gamma_1$ has the
cyclic Rokhlin property and $A_1$ has the unique tracial state,
that $B_1$ has the tracial cyclic Rokhlin property. It follows
that $A_2$ has a unique tracial state. Since $\af_j\otimes {\rm
id}_{\cal U}$ commutes with $\gamma_1,$ it gives an automorphism
on $B_1.$ We denote it by $\af_{j,1}.$ Define
$\gamma_2=\af_{2,1}\otimes \sigma.$ Thus $\gamma_2$ is an
automorphism on $A_2$ satisfying the cyclic Rokhlin property. Let
$G_2\subset K_0(A_2)$ generated by projections with the form
$1_{B_1}\otimes e,$ where $e\in {\cal U}$ is a projection. Then
$$
{\gamma_2}_{*0}|_{G_2}={\rm id}_{G_2}.
$$
Since $A_1$ has a unique tracial state and $\gamma_2$ has the
cyclic Rokhlin property, by applying \ref{PP1}, $B_1$ has a unique
tracial state. We also have that $\rho_{A_2}(G_2)$ is dense in
$\R.$ It follows from \ref{gpots} that $A_2\times_{\gamma_2}\Z$
has tracial rank zero and satisfies the UCT. Moreover, by
(\ref{F1}, we have the following embedding
\beq\label{F2}
(A\rtimes_{\af_1}\Z)\rtimes_{\af_2}\Z \to (B_1\otimes {\cal
U})\rtimes_{\gamma_2}\Z=A_2\rtimes_{\gamma_2}\Z.
\eneq
So, by \ref{2T2}, if $d=2,$ we obtain an embedding from
$A\rtimes_{\Lambda}\Z^d$ into a unital simple AF-algebra.

If $d>2,$ put $B_2=A_2\rtimes_{\gamma_2}\Z.$ We have shown that
$A_2$ has a unique tracial state. Since $\gamma_2$ has the cyclic
Rokhlin property, by \ref{PP1}, $B_2$ has a unique tracial state.
It follows that $A_3=B_2\otimes {\cal U}$ has a unique tracial
state. Note that $\af_{j,1}\otimes {\rm id}_{\cal U}$ is an
automorphism on $A_2.$ Since $\af_{j,1}\otimes{\rm id}_{\cal U}$
commutes with $\gamma_2,$ it gives an automorphism $\af_{j,2}$ on
$B_2.$ Define $\gamma_3=\af_{3,2}\otimes \sigma.$ It is an
automorphism on $A_3$ satisfying the cyclic Rokhlin property. We
then apply the same argument above to conclude that
$A_3\rtimes_{\gamma_3}\Z$ has tracial rank zero and satisfies the
UCT. Moreover, there is an embedding:
\beq\label{F3}
((A\rtimes_{\af_1}\Z)\rtimes_{\af_2}\Z)\rtimes_{\af_3}\Z \to
(B_2\otimes {\cal U})\rtimes_{\gamma_3}\Z
\eneq
In particular, \ref{2T2} implies that the theorem hold for $d=3.$
The theorem follows by applying the same argument and induction.

\end{proof}

\begin{cor}\label{TLc}
Let $A$ be a unital separable simple amenable \CA\, with tracial
rank zero and with a unique tracial state which satisfies the UCT and let 
$G$ be a finitely generated discrete abelian group.
Suppose that $\Lambda: G\to Aut(A)$ is a \hm. Then
$A\rtimes_{\Lambda} G$ can be embedded into a unital simple
AF-algebra
\end{cor}

\begin{proof}
We combine \ref{TL} with an argument of N. Brown. 
Write $G=\Z^d\oplus G_0,$ where $G_0$ is a finite subgroup of $G.$ 
Thus $G/\Z^d\cong G_0$ is compact.
Put $B=A\rtimes_{\Lambda|_{\Z^d}}\Z^d.$ 
By a theorem of Green (Cor. 2.8 of \cite{Gr}), 
$$
A\otimes C(G_0)\rtimes_{\Lambda\otimes \Gamma} G\cong B\otimes M_k
$$
for some integer $k\ge 1.$ 
Now \ref{TL} asserts that $B\otimes M_k$ can be embedded into a unital simple AF-algebra.
The corollary follows from the fact that there is a natural embedding
$A\rtimes_{\Lambda} G\to A\otimes C(G_0)\rtimes_{\Lambda\otimes \Gamma} G$ 
(see also Remark 11.10 of \cite{Bn2}).

\end{proof}

Now we consider crossed products $C(X)\rtimes_{\Lambda}\Z^d$ for
$d>1.$ It is easy to see that if $C(X)\rtimes_{\Lambda}\Z^d$ can
be embedded into a unital simple AF-algebra then there is a
strictly positive $\Lambda$-invariant Borel probability measure.
The following shows that, if in addition, $(X, \af_1)$ is minimal
and unique ergodic for a generator $\af_1$ of $\Lambda,$ then the
converse also holds.

 \begin{thm}\label{LM1}
 Let $X$ be a compact metric space and let $(X,\Lambda)$ be a $\Z^d$ system.
 Suppose
 that there exists a $\Lambda$-invariant strictly positive Borel probability measure.
 Suppose that, in addition, there is a generator $\af_1$ of $\Lambda$ for which
 $(X,\af_1)$ is minimal and unique ergodic.
 Then $C(X)\rtimes_{\Lambda}\Z^d$ can be embedded into a unital simple AF-algebra
  \end{thm}
\begin{proof}
Suppose that $X$ has infinitely many points.
 Let $Y$ be  the
Cantor set. By \ref{LembC} and \ref{C03}, there is a covariant
injective \hm\, $h_1:C(X)\rtimes_{\Lambda}\Z^d\to
C(Y)\rtimes_{\overline{\Lambda}}\Z^d$ such that $(Y,\bt_1)$ is
minimal and unique ergodic.

Therefore, it suffices to prove the theorem in the case that $X$
is the Cantor set.

It follows from 2.1 of \cite{Put}  that $C(X)\rtimes_{\af_1}\Z$ is
a unital simple $A\T$-algebra of real rank zero. Let
$A=C(X)\rtimes_{\af_1}\Z.$ Then $\af_2, \af_3,...,\af_d$ gives a
$\Z^{d-1}$ action $\Lambda'$ on $A.$  By the assumption that $(X,
\af_1)$ is unique ergodic, $A$ has a unique tracial state. It
follows from \ref{TL} that there is an embedding
$$
j:A\rtimes_{\Lambda'}\Z^{d-1}\to C
$$
for some unital simple AF-algebra $C.$  But
$$
A\rtimes_{\Lambda'}\Z^{d-1}\cong C(X)\rtimes_{\Lambda}\Z^d.
$$
Thus we obtain an embedding
$$
C(X)\rtimes_{\Lambda}\Z^d\cong A\rtimes_{\Lambda'}\Z^{d-1}\to C.
$$
\end{proof}

\end{document}